\documentclass{article}

\usepackage{amsmath}
\usepackage{amsthm}
\usepackage{amsfonts}
\usepackage{amssymb}
\usepackage{amscd}
\usepackage{diagrams}

\newtheorem{thm}{Theorem}[section]
\newtheorem{defi}[thm]{Definition}  
\newtheorem{prop}[thm]{Proposition}
\newtheorem{cor}[thm]{Corollary}   
\newtheorem{lem}[thm]{Lemma}
\newtheorem{ex}[thm]{Example}
\newtheorem{rem}[thm]{Remark}

\newcommand{\Q}{\mathbb Q}
\newcommand{\R}{\mathbb R}
\newcommand{\N}{\mathbb N}
\newcommand{\C}{\mathbb C}

\newcommand{\cc}{\mathcal C}

\newcommand{\CG}{\mathcal CG}
\newcommand{\ot}{\otimes}
\newcommand{\s}{\Sigma}

\newcommand{\hook}{\hookrightarrow}
\newcommand{\Map}{\operatorname{Map}}
\newcommand{\Bl}{\operatorname{Bl}}
\newcommand{\de}{\partial}
\newcommand{\ep}{\varepsilon}

\newarrow{Equal} =====
\newarrow{Hook} C--->

\author{Paolo Salvatore \thanks{Supported by the Marie-Curie Research Training Grant ERBFMBICT961611}}
\title{Configuration spaces with summable labels}

\begin{document}

\maketitle

\begin{abstract}
Let $M$ be an $n$-manifold, and let $A$ be a space with a partial sum behaving
as an $n$-fold loop sum.
We define the
space $C(M;A)$ of configurations in $M$ with summable labels 
in $A$ via operad theory. Some examples are symmetric products,
labelled configuration spaces, and
spaces of rational curves.
We show that
$C(I^n,\de I^n;A)$ is an $n$-fold classifying space
of $C(I^n;A)$, and for $n=1$ it is homeomorphic to the classifying
space by Stasheff.
If $M$ is compact, parallelizable, and $A$ is path connected,
then $C(M;A)$ is homotopic to the mapping space
$Map(M,C(I^n,\de I^n;A))$.

\end{abstract}

\section*{Introduction}

The interest in labelled configuration spaces in homotopy theory
dates back to the seventies. May \cite{May} and Segal \cite{Segal2} showed
that  the `electric field map' 
$C(\R^n;X) \to \Omega^n \Sigma^n (X)$ is a weak homotopy equivalence if
$X$ is path connected, and 
in general is the group completion.
Segal showed later \cite{Segal} that the inclusion
$Rat_*(S^2) \hook \Omega^2S^2$
of the space of
based rational selfmaps 
of the spheres into all based selfmaps is the group completion.
He used the identification of $Rat_*(S^2)$ with a space
of configurations in $\C$ with partially summable
labels in $\N \vee \N$, by counting zeros and roots multiplicities.
Guest has recently extended his framework in \cite{Guest}
to the space of based rational curves on projective toric varieties.
Labelled configuration spaces on manifolds have been studied
by B\"odigheimer in
\cite{Bodig}, where the labels are in a based space, and by Kallel in
\cite{Kallel}, where the summable labels belong to a discrete partial abelian monoid.
In both cases the authors have theorems of equivalence between
configuration
and mapping spaces.  
We define configuration spaces on a manifold $M$ with labels in $A$,
where $A$ need not to be abelian. It is sufficient that $A$ has a partial
sum that is homotopy commutative up to level $dim(M)$. The definition
is not trivial and involves tensor products over the Fulton-MacPherson
operad. A substantial part of the paper introduces the necessary tools.
We generalize
the results listed above
to the non-abelian setting, and construct a geometric
$n$-fold delooping in one step.
\
Here is a plan of the paper:

In the first section we define the preliminary notion of a partial
algebra over an operad and its completion.
In the second section we introduce
the Fulton-MacPherson operad $F_n$.
We call an algebra over $F_n$ an $n$-monoid.
A $1$-monoid is exactly an $A_\infty$-space \cite{Kontsevich}. 
In the third section  we describe the homotopical algebra
of topological operads and their algebras.
The main results characterize the homotopy type of $F_n$.

\

{\bf Corollary 3.8.} {\em The unbased Fulton-MacPherson operad $\tilde{F_n}$ is cofibrant.}

\

{\bf Proposition 3.9.} {\em The operad of little $n$-cubes is weakly
equivalent to $F_n$.}

\

This implies that the structure of $n$-monoid is invariant under based
homotopy equivalences, and any connected $n$-monoid has the weak homotopy type of a
$n$-fold loop space. 
In the fourth section
we recall from \cite{Markl} that a partial compactification $C(M)$ of the 
ordered configuration space on an open parallelizable 
$n$-manifold $M$ is a right module over $F_n$.
We define the configuration space 
$C(M;A)$ on $M$
with summable labels in a partial $n$-monoid $A$,
by tensoring $C(M)$ and $A$ over the
operad $F_n$.

The definition of $C(M;A)$ is extended to a general open $n$-manifold $M$
when $A$ is framed, in the sense that $A$ has
a suitable $GL(n)$-action. In the fifth section we define $C(M,N;A)$ for
a relative manifold $(M,N)$ 
by ignoring the particles in $N$. For $n=1$ we obtain the well known
construction by Stasheff:

\

{\bf Proposition 5.11.} {\em If $A$ is a 1-monoid, then
$C(I,\de I;A)$ is homeomorphic to
the classifying space $B(A)$ by Stasheff.}

\

The $n$-monoid completion of a partial $n$-monoid $A$ 
is $C(I^n;A)$, up to homotopy.
We obtain its $n$-fold delooping in one step.

\

{\bf Theorem 6.3.} {\em  If $A$ is framed, then
the group completion of $C(I^n;A)$ is $\Omega^n C(I^n,\de I^n;A)$.}

\

Finally we characterize configuration spaces on manifolds under some
conditions.  

\

{\bf Theorem 6.6.} {\em If $M$ is a compact closed parallelizable $n$-manifold 
and $A$ is a path connected partial framed $n$-monoid,
then there is a weak equivalence $C(M;A) \simeq Map(M;B_n(A))$.}

\

As corollary we obtain a model for the free loop space on a suspension
built out of cyclohedra. This answers a question by Stasheff \cite{operadchik}.

\

{\bf Corollary 6.7.} {\em If $X$ is path connected and well pointed,
then there is a weak homotopy
equivalence $C(S^1;X) \simeq Map(S^1,\s X)$.}

\

The author is grateful to M. Markl and J. Stasheff for many valuable
suggestions.

\section{Partial modules over operads}

Let $\cc$ be a symmetric closed monoidal category, with
tensor product $\otimes$ and unit element $e$. We assume that $\cc$ has 
small limits and colimits.

\begin{defi}
A $\s$-object $X$ in $\cc$ is a collection of objects $X(n)$, for 
$n \in \N$, such that $X(n)$ is equipped with an action of the
symmetric group $\s_n$.
\end{defi}

The category of $\s$-objects in $\cc$ will be denoted by $\s \cc$.
We observe as in \cite{GJ} that 
$\s \cc$ is a monoidal category as follows:
given two objects $A$ and $B$, their tensor product $A \otimes B$ is
defined by
$$(A \otimes B)(n)  = \coprod_{k=0}^{\infty}A(k) \ot_{\s_k} 
(\coprod_{\pi \in \Map(n,k)} \bigotimes_{i=1}^{k} B_{\pi^{-1}(i) }\, )  ,$$
where $Map(n,k)$ is the set of maps from $\{1,\dots,n\}$ to
$\{1, \dots, k \}$.
Here each element $B_S$, where $S$ is a set of numbers,
is identified to $B_{\#S}$ by the order preserving bijection, and
the action of $\Sigma_n$ is given accordingly.
There is a natural embedding functor $j:\cc \hookrightarrow \s \cc$
considering an object $X$ as a $\s$-object concentrated in degree 0, so that   
 
\begin{equation*}
j(X)(n)=
\begin{cases}
X& \text{if $n=0$}, \\
\emptyset& \text{if $n \neq 0$}.
\end{cases}
\end{equation*}

with the trivial actions of the symmetric groups. Here $\emptyset$ denotes the
initial object of $\cc$.
The functor $j$ is left adjoint to the forgetful functor
$B \mapsto B(0)$ from $\s \cc$ to $\cc$.
More generally we have an embedding functor $j_n: \s_n  \cc \to \s \cc$,
that is left adjoint to the forgetful functor $B \mapsto B(n)$.

The unit  element $\iota$ of $\s \cc$ is defined by
\begin{equation*}
\iota(n)=
\begin{cases}
e& \text{if $n=1$}, \\
\emptyset& \text{if $n \neq 1$}.
\end{cases}
\end{equation*}

\begin{defi} {\rm \cite{GJ}}
An operad in $\cc$ is a monoid in the monoidal category $\s \cc$.
We denote the category of operads in $\cc$ by $\mathcal{OP}(\cc)$.
\end{defi}

This means that an operad $(F,\mu, \eta)$ is a $\s$-object $F$
together with composition morphism $\mu: F \otimes F \to F$ and unit
morphism $\eta: \iota \to F$, such that the associativity property 
$\mu (\mu \otimes F) = \mu (F \otimes \mu): F \otimes F \otimes F \to F$  and the unit
property $\mu (F \otimes \mu) = \mu (\mu \otimes F) = id_F :F \to F$ hold.
Note that the functor $F \otimes \_\;$ forms a triple.

\begin{ex}
The category $\mathcal{CH}_R$ of non negatively graded chain complexes over a commutative ring $R$
is monoidal by the tensor product. The operads in  $\mathcal{CH}_R$
are called
differential graded operads over $R$.
\end{ex}

\begin{defi} {\rm\cite{Markl}}
Given an operad $F$ in $\cc$,
a left $F$-module $A$ is
a $\s$-object $A$, with a morphism 
$\rho: F \otimes A \to A$ of $\s$-objects such that

$ \rho (F \otimes \rho) = \rho (\mu \otimes A) : F \otimes F \otimes A \to A$
and $\rho (\eta \otimes A) = id_A:A \to A$.
\end{defi}
In other words $A$ is an algebra over the triple $F \otimes \_\;$.
Dually we define the notion of right 
$F$-module. 
We denote the category of left $F$-modules by $Mod_F$, and the
category of right $F$-modules by $_F Mod$.
\begin{defi}
If $F$ and $G$ are operads in $\cc$, then a $F$-$G$-bimodule $A$
is a left $F$-module and a right $G$-module such that the left $F$-module
structure map is a right $G$-homomorphism.
\end{defi}
\begin{defi}
An algebra $X$ over an operad $F$, or $F$-algebra, is an object $X$ of $\cc$ together with
a left $F$-module structure on $j(X)$.
\end{defi}
We denote the category of $F$-algebras by $Alg_F$.
Moreover we will denote by $F(Y)$ the free $F$-algebra
generated by the object $Y$ in $\cc$. This object is defined by 
$j(F(Y)) = F \ot j(Y)$.

\begin{defi} 
A partial left $F$-module $A$ is a $\s$-object $A$ in $\cc$ together with
a monomorphism $i:Comp \hookrightarrow F \otimes A$ in $\s \cc$
and a composition map $\rho: Comp \to A$ such that
\begin{enumerate}
\item The unit $\eta \otimes A: A \to F \otimes A$ factors uniquely
through $\tilde{\eta}: A \to Comp$ and 
the composition $\rho (\tilde{\eta})= id_A$ is the identity.

\item The pullbacks in $\s\cc$ of $i:Comp \hookrightarrow F \otimes A$ along the two
maps
$(\mu \otimes A)(F \otimes i) : F \otimes Comp \to F \otimes A$
and 
$F \otimes \rho : F \otimes Comp \to F \otimes A$
coincide. Moreover the compositions of the two pullback maps 
with $\rho: Comp \to A$ coincide.
\end{enumerate}
\end{defi}

Partial right $F$-modules are defined dually.
A partial $F$-algebra is an object $X$ of $\cc$ such that 
$j(X)$ is a partial left $F$-module.
We denote the categories of partial left $F$-modules, right $F$-modules
and $F$-algebras respectively by
$PartMod_F , \; _F PartMod$, and  $PartAlg_F$.

A morphism $g:(A,Comp_A) \to (B,Comp_B)$ of  partial left $F$-modules is a morphism in 
$\s \cc$ such that $(F \otimes g) i_A:Comp_A \to F \otimes B$ factors
through
$\tilde{g}:Comp_A \to Comp_B \;$, and
$g \rho_A  = \rho_B \tilde{g}$. 

\

We exhibit a functor from partial to total left modules that is left
adjoint to the forgetful functor. The analogous construction for 
right modules is exactly dual.

If $A$ is a partial $F$-left module, then define $\hat{A}$ by
 the coequalizer in the category $\s\cc$
$$ F \otimes Comp
 \pile{\rTo^{(\mu  \otimes A)(F \otimes i)} \\
\rTo_{F \otimes \rho} }
F \otimes A \rDotsto \hat{A} \; .$$

\begin{prop}
There is a left $F$-module structure on $\hat{A}$.
\end{prop}

\begin{proof}
The proof is modelled on Lemma 1.15  in \cite{GJ}.
The coequalizer above is reflexive because the input arrows admit the
common section
$F \otimes \tilde{\eta}: F \otimes A \rightarrow F \otimes Comp$. 
Now $\hat{A}$ admits the
structure of left $F$-module, because by
Lemma 2.3.8 in \cite{Rezk} $F \ot \_$ preserves reflexive coequalizers.
Moreover $\hat{A}$
is the coequalizer of the pair
above in the category of left $F$-modules.
\end{proof}

\begin{prop}
The completion $A \mapsto \hat{A}$ induces a functor that is 
left adjoint to
the forgetful functor $U: Mod_F \to PartMod_F$.
\end{prop}

\begin{defi}
For any right $F$-module  $C$ with structure map
$\sigma: C \otimes F \to C$ and a partial left $F$-module $(A,Comp,i)$
 we define the tensor product
$C \otimes_F A$ as coequalizer in $\s\cc$

$$ C \otimes Comp_A
\pile{\rTo^{(\sigma \otimes A)(C \otimes i)}   \\ \rTo_{C \otimes \rho}}
C \otimes A \rDotsto C \otimes_F A \, .$$

\end{defi}
Dually we define the tensor product of a partial right $F$-module and a left $F$-module.

\begin{prop} \label{tensor}
Given a partial right $F$-module $A$, an $F$-$G$-bimodule $B$, and
a partial left $G$-module $C$, there are natural isomorphisms
$$(A \ot_F B) \ot_G C \cong A \ot_F (B \ot_G C).$$
\end{prop}

The isomorphism holds because the tensor product is a left adjoint
and preserves colimits.

\section{The Fulton-MacPherson operad}
The category $\mathcal{CG}$ of compactly generated weak Hausdorff topological spaces
is a closed monoidal category with all limits and colimits, hence it
satisfies the assumptions of the previous section. We note however that
in general the forgetful functor to the category of sets does not preserve
colimits. 
Operads and modules in $\mathcal{CG}$ shall be called simply topological
operads and topological modules.

\

The key topological operads in this paper are the Fulton-MacPherson operads,
that are suitable cofibrant versions of the little cubes operads.
They were introduced in \cite{GJ}. We recall their definition.
Consider the differential-geometric blow-up of $(\R^n)^k$ along the
small diagonal $\Delta = \{ x_1,\dots,x_k \: | \: x_1 = \dots = x_k \}$.
The blow-up is explicitly obtained if we replace the diagonal by its
normal sphere bundle. The fiber of the trivial normal bundle at the origin is 
$F= \{ y_1, \dots, y_k  \:|\: \sum_{i=1}^{i=k}y_i = 0   \}$ 
and the sphere bundle 
$PF = (F-0)/(\R^{+})$ can be seen as the space of closed half-lines in $F$.
Then the blow-up is
$$\Bl_{\Delta}((\R^n)^k) = \{ (x,y) \in (\R^n)^k \times PF \, |\,
x-\pi_\Delta(x) \in  y \}, $$
where the orthogonal projection is $\pi_\Delta(x_1,\dots,x_k) =
(\sum_{i=1}^{i=k} x_i /k, \dots ,\sum_{i=1}^{i=k} x_i /k).$
For any set $S \subseteq \{1,\dots,k\}$ let us denote by
$\Bl_\Delta((\R^n)^S)$  the blow-up of $(\R^n)^S$ along its small diagonal.

\

Let $C_k^0(\R^n) \subset Map(\{1,\dots,k\},\R^n)$ be the space of ordered pairwise
distinct $k$-tuples in $\R^n$. There is a natural right $\s_k$-action on
this space, and we consider it as left $\s_k$-space by the opposite action.
 As  $C_k^0(\R^n)$ does not intersect any
diagonal, there is a natural embedding
$$j: C_k^0(\R^n) \to \prod_{S \subseteq \{1,\dots,n\}, \,\#S \geq 2}
 \Bl_\Delta ((\R^n)^S).$$

\begin{defi}
The Fulton-MacPherson configuration space
$C_k(\R^n)$ is the
closure of the image of $j$.
\end{defi}
We note that $GL(n)$ acts diagonally on each blowup,
$j$ is a $GL(n)$-equivariant map and therefore
$C_k(\R^n)$ is a $GL(n)$-space.
\

In a similar way we define the Fulton-MacPherson
configuration space $C_k(M)$ of a
smooth open manifold $M$.
In this case one builds the differential-geometric blowups
of $M^k$ along the diagonal $\Delta_M$ by gluing together
$M^k - \Delta_M$ and the normal sphere bundle via a tubular
neighbourhood of $\Delta_M$ in $M^k$.
It turns out that $C_k(M)$ is a manifold with corners $\s_k$-equivariantly homotopy
equivalent to its interior, the ordinary configuration
space $C_k^0(M)$ of ordered pairwise distinct $k$-tuples in $M$.
                                            
There is a blow-down map $b:C_k(M) \to M^k$ such that the composite
$C_k^0(M) \rHook^j C_k(M) \rTo^b M^k$
is the inclusion. We will say that the blow-down map gives the {\em macroscopic
 locations} of the
particles.

There is a characterization of the Fulton-MacPherson configuration
space by means of trees due to
Kontsevich.
For us a {\em tree} is an oriented finite connected 
graph with no cycles such that each vertex has exactly
one outcoming edge. An {\em ordered tree} is a tree together with 
an ordering of the incoming edges of each vertex.
The ordering is equivalent to the assignation of a planar embedding.
The only edge with no end vertex is the
{\em root}, the edges with no initial vertex are the {\em twigs},
and all other edges are {\em internal}. A tree on a set $I$ is a tree
together with a bijection from the set of its twigs to $I$.
The {\em valence} of a vertex is the number of incoming edges.
Let $G(n)$ be the group of affine transformations of $\R^n$
generated by translations and positive
dilatations.

\begin{prop} {\rm \cite{Kontsevich}} \label{Kont}
Let $M$ be an open manifold. Then each element
in $C_k(M)$ is uniquely determined by:
\begin{enumerate}
\item Distinct macroscopic locations
$P_1,\dots,P_l \in M$, with $1 \leq l \leq k$.
\item For each $1 \leq i \leq l$ a tree $T_i$ with $f_i$ twigs, so that
$\sum_{i=1}^l f_i = k$, and for each vertex in $T_i$ of valence $m$
an
element in  $C^0_m (\tau_{P_i}(M))/G(n)$, where $\tau_{P_i}(M)$ is the tangent
plane at $P_i$.
\item A global ordering of the $k$ twigs of the trees.
\end{enumerate}
\end{prop}

\begin{defi}
If $b:C_k(\R^n) \to (\R^n)^k$ is the blowdown map, then 
the Fulton-MacPherson
space is  $F_n(k) = b^{-1}(\{0\}^k)$.
\end{defi}
This space contains all configurations macroscopically located
at the origin.
\begin{prop} \rm{\cite{GJ}}
The space $F_n(k)$ is a manifold with corners, and it is a
compactification of $C^0_k(\R^n)/G(n)$.
\end{prop}
The {\em faces} of $F_n(k)$ are indexed by trees on $\{1,\dots,k\}$,
and the codimension of a face is equal to the number of internal edges
of the indexing tree.

\begin{prop} {\rm \cite{Markl}}
The spaces $F_n(k)\; k \geq 0$ assemble to form a topological operad.
Moreover $F_n(k)$ is a $\s_k$-equivariant deformation retract
of $C_k(\R^n)$.
\end{prop}

The composition law is easily described in terms of trees:
each element in $F_n$ is described by a single tree by Proposition \ref{Kont}.
If $a \in F_n(k) $ and $b_j \in F_n(i_j)$ for $j=1, \dots, k$ 
then $a \circ (b_1,\dots,b_k) \in F_n(i_1 + \cdots + i_k)$ corresponds
to the tree obtained by merging the $j$-th twig of the tree of $a$ with
the root of the tree of $b_j$ for $j = 1, \dots,k$, and assigning the new
twigs the induced order. This operation on trees will be called
{\em grafting}.
Note that $F_n(1)$ is a point, the unit $\iota$ of the operad, and
is represented by the trivial tree.
We assume that $F_n(0)$ is a point, the empty configuration.
We stress the fact that Getzler and Jones in \cite{GJ} consider the unpointed version
$\tilde{F}_n$ such that $\tilde{F}_n(k)=F_n(k)$ for $k>0$ and $\tilde{F}_n(0)=\emptyset$.
Their paper focuses on the differential graded operad 
$e_n= H_*(\tilde{F}_n;\Q)$, the rational homology
of $\tilde{F}_n$. They denote $H_*(F_n;\Q)$ by $e_n^+$.
The deformation retraction $r:C_k(\R^n) \times I \to C_k(\R^n)$ such that
$r(C_k(\R^n) \times \{1\}) = F_n(k)$ 
is defined for $t \neq 0$ and $x \in C_k^0(\R^n)$ by 
$r(x,t) =  xt $.
\

\begin{defi}
We call an algebra over $F_n$ an $n$-monoid, and 
an algebra over $\tilde{F}_n$ an $n$-semigroup.
\end{defi}

\begin{ex} \rm{\cite{Kontsevich}}
The 1-monoids are the $A_\infty$-spaces.
\end{ex}

In fact $F_1(i) = K_i \times \s_i$, where $K_i$ denotes
the associahedron by Stasheff \cite{operadchik}, so $F_1$ is the symmetric
operad generated by the non-symmetric Stasheff
operad $K$. But an $A_\infty$ space is
by definition an algebra over $K$.

\section{Homotopical algebra and the little discs}

We describe the
closed model category structure of
the categories of topological operads and their 
algebras.

\begin{defi} {\rm \cite{DKH}}
A cofibrantly generated model category
is a closed
model category {\rm \cite{Quillen}}, together
with a set $I$
of {\em generating cofibrations}, and a set $J$ of
{\em generating trivial cofibrations}, so that
the fibrations and the trivial fibrations are respectively
the maps satisfying the right lifting property with respect to the maps
in $J$ and $I$.
\end{defi}

Consider the free operad functor
$\mathbb{T}:\s(\mathcal{CG}) \to \mathcal{OP}(\mathcal{CG})$,
left adjoint to 
the forgetful functor $\mathbb{U}: \mathcal{OP}(\mathcal{CG}) \to
\s(\mathcal{CG})$.
Let $\mathcal{S}_n$ be the family of subgroups of $\s_n$.
The simplicial version of the following
proposition is 3.2.11 in \cite{Rezk}.

\begin{prop} {\rm \cite{Salvatore}}  \label{pacco}
The category of topological operads is a cofibrantly generated model
category, with the following structure:

\begin{enumerate}

\item The set of generating cofibrations
is
$I = \{ \mathbb{T}(
\de I^i \times H \setminus \s_n \hook I^i \times H \setminus \s_n )
   \; | \;  i,n \in \N, \, H \in \mathcal{S}_n \,  \};$

\item The set of generating
trivial cofibrations is

$J = \{ \mathbb{T} ((I^{i-1} \times \{0\})  \times H \setminus \s_n \hook I^i \times H\setminus
\s_n ) \; | \; i,n \in \N, \, H \in \mathcal{S}_n \,  \};$

\item A morphism $f$ is respectively a weak equivalence or a fibration if 
for any $n \in \N$ and $H \in \mathcal{S}_n$
the restriction $f_n^H$ of $f_n$ to the $H$-invariant subspaces
is respectively a weak homotopy equivalence or a Serre fibration.

\end{enumerate}
\end{prop}

There is a functorial cofibrant resolution for operads, introduced
in \cite{Boardman}.
Let $A$ be a topological operad.
Let $M_k$ be the set of isomorphism classes of {\em ordered} trees
on $\{1,\dots,k\}$, and for each tree $t$ let $V(t)$ be the set
of its vertices and $E(t)$ the set of its internal edges.
For each vertex $x \in V(t)$ let $|x|$ be its valence.

\begin{defi}
The space
of ordered
trees on $\{1,\dots,k\}$ with vertices
labelled by elements of $A$, and with internal edges 
labelled by real
numbers in $[0,1]$ is
$$T_k(A) = \coprod_{t \in M_k} ( \prod_{x \in V(t)} A(|x|) \times [0,1]^ {\# E(t)} )\,.$$
\end{defi}

Let $T_t$ be the summand indexed by a tree $t \in M_k$.
For each internal edge $e \in E(t)$ the operad composition
induces a map $\de_e: T_t \to T_{t-e}$,
where $t-e$ is obtained from $t$
by collapsing $e$ to a vertex.
If $e$ goes from $x$ to $y$, 
$|y| = n$, and $e$ is the $i-th$ incoming edge of $y$, then
$\de_e(x)$ is the multiplication of
the composition $\theta_i: A(|x|) \times A(|y|) \longrightarrow
A(|x|+|y|-1)$ by the identity maps of the vertices in $V(t)-\{x \cup y\}$.

\begin{defi}
The space $WA(k)$ is the quotient of 
$T_k(A)$ under the following relations:

\begin{enumerate}

\item  Suppose that $t \in T_k(A)$, $v$ is a vertex of $t$ of valence $n$
labelled by $\alpha \in A(n)$, the subtrees
stemming from $v$ are $t_1 < \dots <t_n$, and $\sigma \in \s_n$.
Then $t$ is equivalent to the element obtained from $t$ by replacing
$\alpha$ by $\sigma^{-1} \alpha$ and by permuting the order of the subtrees to
$t_{\sigma_1}< \dots < t_{\sigma_n}$.

\item \label{relation} If $t \in T_k(A)$ has an edge $e$ of length 0,
then $t$ is equivalent to the labelled tree obtained
by collapsing $e$ to a vertex, and composing the labels of its
vertices.

\item \label{unita} If $t \in T_k(A)$ has a vertex $w$ of valence 1 labelled
by the unit $\iota \in A(1)$
of the operad $A$, then $t$ is equivalent to the labelled tree obtained
by removing $w$. If $w$ is between two internal edges of lengths $s$ and $t$,
then we assign length $s + t - st$ to the merged edge.
\end{enumerate}
\end{defi}

There is an action of $\s_k$ on $WA(k)$ induced by permuting the
labelling of the twigs of elements in $T_k(A)$. 

\begin{prop} {\rm \cite{Boardman}}
There is an operad structure on $WA$, defined
by grafting trees, and by assigning length
1 to the
new internal edges. A natural ordering of the twigs of the composite
is induced.
The trivial tree consisting of an edge with no vertices is the identity of $WA$.
\end{prop}

For us a cofibration of topological spaces is the retract of
a generalized CW-inclusion \cite{DKH}. 
We say that a pointed space $(X,x_0)$ is well-pointed if the inclusion
$\{x_0\} \hook X$ is a cofibration.
The following proposition is essentially proved in \cite{Boardman}.
\begin{prop}
Let $A$ be a topological operad such that $(A(1),\iota)$ is well-pointed
and each $A(n)$ is a cofibrant space.
Then $WA$ is a cofibrant resolution of $A$.
\end{prop}

\begin{prop} \label{rosso}
There is an isomorphism of topological operads $W(\tilde{F}_n)  \cong \tilde{F}_n$.
\end{prop}

\begin{proof}

We observe that $W (\tilde{F}_n)(i)$ is obtained by gluing together
for each {\em face} $S$ of $F_n(i)$ of codimension $d$ a copy
of
$S  \times [0,1]^d$.
This is true because the codimension of a face of the manifold with corners
$F_n(i)$ is equal to the number of internal edges of the associated tree.
Then
$W (\tilde{F}_n)(i)$ admits
the structure of manifold with corners diffeomorphic to 
$F_n(i)$, and the composition maps of both operads are
described by grafting of trees.
\end{proof}
\begin{cor}
The operad $\tilde{F}_n$ is cofibrant.
\end{cor}

Let $D_n$ be the operad of little $n$-discs.
The space $D_n(k)$ is the space of $k$-tuples of direction preserving
affine selfembeddings
of the unit $n$-disc with pairwise disjointed images.
The operad structure is defined by multicomposition of the embeddings.
There is a sequence of $\s_k$-equivariant
homotopy equivalences
$D_n(k) \to C_k^0(I^n) \hook C_k^0(\R^n) \hook C_k(\R^n) \to F_n(k).$
The first map sends the little discs to
their centers, the last is a deformation retraction.
The inclusion $C_k^0(I^n) \hook C_k^0(\R^n)$ is a $\s_k$-equivariant
homotopy equivalence because the inclusion $I^n \hook \R^n$ is isotopic
to a homeomorphism.
The image of the composite $r_k$ is
the interior of the manifold with
corners $F_n(k)$. It follows that the $\s$-map
$r=\{r_k\}$ is not an operad map
because all elements in the boundary of $F_n(k)$ are composite.

\begin{prop} \label{cina}
The operad $D_n$ is weakly equivalent to $F_n$.
\end{prop}

\begin{proof}
We build an extension $R: WD_n \to F_n$ of $r:D_n \to F_n$ that is
a map of operads and a weak equivalence.

An element $a$ in $WD_n(k)$ is represented by a labelled tree $\tau \in
T_k(D_n)$ on
$\{1,\dots,k\}$. If the $i$-tuple $(f_1^v,\dots,f_i^v)$
labels a vertex $v$ of valence $i$, then for each $j$ we associate the
embedding $f_j^v$ to the $j$-th {\em incoming}
 edge $e_j(v)$ of $v$. The equivalence relation defining $WD_n$
preserves this association.
We observe incidentally that the multicomposition of the labels of $\tau$
is the $k$-tuple of embeddings $(g_1,\dots,g_k)$ such that for each $j\,$
$g_j$ is  the composition
of the embeddings associated to the edges along the unique path
from the $j$-th twig to the root.
Suppose that the internal edges of $\tau$
are labelled by numbers in $(0,1)$.
Let $l(e)$ denote the length of an edge $e$ and
if $r \in (0,1]$ let $\delta_r$ be the dilatation of the $n$-disc by  $r$.
Consider the labelled tree $\tau'$  obtained from $\tau$ by replacing
for each vertex $v$ and for each $j=1,\dots,|v|$
 the embedding $f_j^v$ by the rescaling $f_j^v  \circ \delta_{l(e_j(v))}$.
Let $b$ be the multicomposition of the labels of $\tau'$, and set
$R_k(a)=r_k(b)$.
We have defined $R_k$ on a dense subspace
of $WD_n(k)$. The map $R_k$ extends 
to $WD_n(k)$ and $R$ is an operad map, because the boundary and the composition
of $F_n$ are described by a limit procedure.
Let $i_k:D_n(k) \to WD_n(k)$ be the inclusion such that $i_k(a)$ is
represented by the tree on $\{1,\dots,k\}$
with a single vertex labelled by $a$. 
The map $R_k:WD_n(k) \to F_n(k)$ is a $\s_k$-equivariant homotopy equivalence
for each $k$ because $i_k$ is such \cite{Boardman} and $R_k i_k = r_k$. 
In particular $R$ is a weak equivalence
of topological operads, and $D_n \simeq F_n$.
\end{proof}

The simplicial analogue of the following proposition is
3.2.5 in \cite{Rezk}.  

\begin{prop} {\rm \cite{Salvatore} \cite{SV}}
Let $F$ be a topological operad. Then the category $Alg_F$ is a
cofibrantly generated model category with the following structure:

\begin{enumerate}

\item The set of generating cofibrations is $I = \{ F ( \de I^i ) \hook
F ( I^i )\, | \, i \in \N  \}$.

\item  The set of trivial generating cofibrations is
$J = \{ F  (I^{i-1} \times \{0\}) \hook F ( I^i )\, | \,  i \in \N  \}$.

\item A $F$-homomorphism is a weak equivalence or a fibration
if it is respectively a weak homotopy equivalence or a Serre fibration.

\end{enumerate}
\end{prop}

Under mild conditions there is a functorial cofibrant resolution of
topological algebras over operads, introduced in \cite{Boardman}.
Let $A$ be a topological operad.
Consider the $\s$-space  $W^+A$ defined similarly as $WA$, except that
relation \ref{unita} is not applied if $w$ is the root vertex.
It turns out that $W^+A$ is an $A-WA$-bimodule, by action of
$A$ on the label
of the root, and by grafting trees representing elements of $WA$.
Let $X$ be an $A$-algebra. It has a $WA$-algebra structure induced by the
projection $\ep: WA \to A$.
We define the $A$-algebra $U_A(X) = W^+A \ot_{WA} X$.
The projection $W^+ A \to WA$,  obtained by extending relation \ref{unita}
to the root vertex, induces an $A$-homomorphism $\pi: U_A(X) \to X$,
that is a deformation retraction, see p. 51 of \cite{Boardman}.

\begin{prop} {\rm \cite{Salvatore}}
If $X$ is a cofibrant space
then $U_A(X)$ is a cofibrant $A$-algebra.
\end{prop}

\begin{defi}
If $A$ is a topological operad, $X$ and $Y$ are $A$-algebras, then
a homotopy $A$-morphism from $X$ to $Y$ is
an $A$-homomorphism from $U_A(X)$ to $Y$.
\end{defi}

\begin{prop}
If $F$ is a topological operad, $X$ is the retract of a generalized $CW$-space,
and $Ho(Alg_F)$ is the homotopy category, then
$Ho(Alg_F)(X,Y) = Alg_F(U_F(X),Y)/\simeq$.
\end{prop}

\begin{proof}
The set of right homotopy classes $[X,Y]$ in the sense of \cite{Quillen} is
the set  of
$F$-homomorphisms from
a cofibrant model of $X$ to a fibrant model of $Y$
 modulo right homotopy.
Now $U_F(X)$ is a cofibrant resolution of $X$, and $Y$ is fibrant
because every object
is such. It is easy to see that the
the right homotopy classes of $F$-homomorphisms from $U_F(X)$ to $Y$ are ordinary
homotopy classes, because $Y^I$ is a path object.
\end{proof}

This result is consistent with the formulation of the homotopy category
of $F$-algebras in \cite{Boardman}.

\begin{prop} 
If $Z$ is an $n$-semigroup and $p:Y \to Z$ is a homotopy equivalence, then
$Y$ has a structure of $n$-semigroup such that $p$ extends to a
homotopy
$\tilde{F}_n$-morphism.
\end{prop}
\begin{proof}
It is sufficient to observe that $W \tilde{F}_n $ is homeomorphic to
$\tilde{F}_n$, and apply the homotopy invariance theorem 8.1 in
\cite{Boardman}.
\end{proof}

\section{Modules and configuration spaces with summable labels}

\begin{prop} {\rm \cite{Markl}}
 For any parallelizable open
manifold $M$ of dimension $n$, the space of configurations $C(M)= \coprod_{k \in \N}C_k(M) $
is a right $F_n$-module.
\end{prop}
\begin{proof}
We choose a trivialization
of the tangent bundle $\tau(M) \cong M \times \R^n$. Then
the composition $C(M) \ot F_n \to C(M)$ is described by grafting of trees
representing elements as in \ref{Kont}.
\end{proof}

Markl in \cite{Markl} gives a similar picture for generic open manifolds by
introducing the framed Fulton-MacPherson operads.

\begin{defi}
Let $G$ be a topological group, and let $F$ be a topological operad
such that $F(i)$ is a $G \times \s_i$-space for each $i$, and
the structure map $\mu$ of $F$ is $G$-equivariant.
The semidirect product
$F \rtimes G$ is the operad defined by
$(F \rtimes G )(i) = F(i) \times G^i$, with structure map

\begin{equation*}
\begin{split}
\tilde{\mu}((x,g_1,\dots,g_k);(x_1,g_1^1,\dots,g_1^{m_1}),
\dots,(x_k,g_k^1,\dots,
g_k^{m_k})) = \\
= ( \mu(x;g_1x_1,\dots,g_k x_k),g_1 g^1_1,\dots ,g_k g_k^{m_k}) \;. \\
\end{split}
\end{equation*}
\end{defi}

\begin{defi}
The framed Fulton-MacPherson operad is the semidirect product
$fF_n = F_n \rtimes GL(n)$.
\end{defi}

\begin{defi} 
Let $M$ be an open $n$-manifold.
The $GL(n)$-bundle of frames on $M$ induces a
$GL(n)^k$-bundle $fC_k(M)$ on $C_k(M)$, acted on by $\s_k$, that we call
the {\rm framed configuration space} of $k$ frames in $M$.
\end{defi}

\begin{prop} {\rm \cite{Markl}}
The $\s$-space $fC(M)$ of framed configurations is
a right module over $fF_n$.
\end{prop}

An element of the framed configuration space
$fC_k(M)$ is uniquely determined by labelled
trees 
as in Proposition \ref{Kont}, and by additional $k$ frames
of the tangent planes associated to the $k$ twigs.
A smooth embedding $i:M \hook N$ of open $n$-manifolds induces
a right $fF_n$-homomorphism $fC(i):fC(M) \hook fC(N)$.
\begin{rem}
If $M$ is a Riemannian $n$-manifold, then
we can define for each $k$
a $O(n)^k$-bundle $f^O C_k(M)$ over $C_k(M)$,
so that $f^O C(M)$ is a right $F_n \rtimes O(n)$-module.
If $M$ is oriented then we can define a $SO(n)^k$-bundle $f^{SO}C_k(M)$
on $C_k(M)$ so that $f^{SO} C(M)$ is a right $F_n \rtimes SO(n)$-module.

\end{rem}

\begin{defi}
We call an algebra over $fF_n$ a framed $n$-monoid.
\end{defi}
Hence a framed $n$-monoid
is an $n$-monoid equipped with an action of $GL(n)$, that is
compatible with the $n$-monoid structure map.

\begin{defi}
A partial framed $n$-monoid is a partial $n$-monoid with an action of
$GL(n)$, such that $GL(n)$ preserves $Comp$ and respects
the partial composition.
\end{defi}

\begin{defi}
Let $fD_n(k)$ be the space of $k$-tuples of
of affine selfembeddings of the unit $n$-disc that preserve
angles and have
pairwise disjointed images.
The multicomposition gives $fD_n$
the structure of an operad, that we call the operad
of framed little $n$-discs.
\end{defi}

\begin{rem}
Consider the iterated loop space $\Omega^n(X,x_0)$ as the space of
maps from the closed unit $n$-disc to $X$, sending the boundary to 
the base point $x_0$. This space is an algebra over $fD_n$. 
\end{rem} 

\begin{prop}
The operad of framed little $n$-discs $f D_n$ is weakly equivalent
to the framed Fulton-MacPherson operad $fF_n$.
\end{prop}
\begin{proof}
We apply the same proof of
Proposition \ref{cina}
to show that $fD_n \simeq F_n \rtimes O(n)$, and conclude by
the homotopy equivalence
$O(n) \hookrightarrow GL(n)$.
\end{proof}

If we restrict to the suboperad $\bar{f} D_n \subset fD_n$ containing
orientation preserving embeddings, then we obtain a weak equivalence
$\bar{f} D_n \simeq F_n \rtimes SO(n)$.

\begin{defi}  \label{previous}
Let $A$ be a partial $n$-monoid, and let $M$ be
an open parallelizable manifold of dimension $n$.
Then the space of 
configurations in $M$ with partially summable labels in $A$ is
$C(M;A) : = C(M) \ot_{F_n} A$.
\end{defi}

An element of $C(M) \ot A = \coprod_k C_k(M) \times_{\s_k}A^k$
consists by \ref{Kont} of a finite set of trees
based at distinct points in $M$, with vertices labelled by $F_n$ and
twigs labelled by $A$. The equivalence relation defining $C(M;A)$
says that if some twigs labelled by $a_1,\dots,a_k$ 
are departing from a vertex labelled by $c \in F_n(k)$ in 
$t \in C(M) \ot A$ and
$\rho(c;a_1,\dots,a_k)$ is defined, then we identify $t$ with the forest
obtained from $t$ by cutting such twigs, and by replacing their vertex
by a twig labelled by $\rho(c;a_1,\dots,a_k)$. Furthermore if the $i$-th
twig departing from a vertex labelled by $c$ in $t$
is labelled by the base point $a_0$ then we identify $t$ to forest obtained
by cutting the twig and by replacing the label $c$ by $s_i(c)$,
where 
$s_i:F_n(k) \to F_n(k-1)$ is the projection induced by forgetting the
$i$-th coordinate.

If $A$ is an $n$-monoid, then by iterated identifications any element
in $C(M;A)$ has a unique representative consisting of a finite
set of trivial trees in $M$, or points, with labels in $A-\{a_0\}$.


We denote by $|\quad|: \CG \to {\bf Set}$ be the forgetful functor.

\begin{prop} \label{pish}      
Suppose that the inclusion $Comp \hook F(A)$ is a cofibration,
and $A$ is well-pointed. Then

\begin{enumerate}
\item $|C(M;A)| = |C(M)| \ot_{|F_n|}|A| \;$;

\item  the space $C(M;A)$ has the weak topology with respect
to the filtration  $C_k(M;A) = Im( \coprod_{i \leq k} C(M)_i
\times_{\s_i} A^i) ,\;$
$k \in \N$.
\end{enumerate}
\end{prop}

\begin{proof}
If $A$ is a proper $n$-monoid, then
we have relative homeomorphisms 
$$(C_k(M),\de C_k(M)) \times_{\s_k} (A,a_0)^k \longrightarrow
(C_k(M;A),C_{k-1}(M;A)    )$$ for $k \geq 1$, and
we conclude by 8.4, 9.2 
and 9.4 in \cite{Steenrod}.
If $A$ is a partial $n$-monoid, then 
we denote by 
$R_i \subset C_i(M) \times_{\s_i} A^i$ the space
 of reducible elements that are equivalent
to an element of some $C_j(M) \times_{\s_j} A^j$ with $j<i$.
For example,
\begin{align*}
R_1 & = M \times \{a_0\} ; \\
R_2 & = (C_2(M) \times_{\s_2} ( A \vee A)) \cup (M \times Comp_2 ) ; \\
R_3  & = (C_3(M) \times_{\s_3}(A\vee A\vee A))  \cup
(C_2(M) \times_{\s_2} (Comp_2 \times A)) \cup (M \times Comp_3) .\\
\end{align*}
We have relative homeomorphisms 
$( C_i(M) \times_{\s_i} A^i, R_i) \to (C_i(M;A) , C_{i-1}(M;A))\, ,$ 
 and we argue similarly.
\end{proof}

\begin{defi}
Suppose that $M$ is an open $n$-dimensional smooth manifold,
and $A$ is a partial framed $n$-monoid.
Then the space of configurations in $M$ with labels in $A$ is 
 $C(M;A):= fC(M) \otimes_{fF_n} A$.
\end{defi}

Note that if $M$ is parallelizable then the definition is consistent with
\ref{previous}. In fact the framed configurations in $M$ are
given by $fC(M) = C(M) \ot_{F_n} fF_n$ and  by \ref{tensor}
$$fC(M) \ot_{fF_n} A = C(M) \ot_{F_n} fF_n \ot_{fF_n} A = C(M) \ot_{F_n} A.$$

\begin{prop}
Let $A$ be a partial framed $n$-monoid with base point $a_0$
such that the inclusions
$Comp \hookrightarrow fF_n (A)$ and 
$\{a_0\} \hook A$ are
cofibrations of $GL(n)$-spaces.
 Let $M$ be an open $n$-manifold.
Then
\begin{enumerate}

\item $|C(M;A)| = |fC(M)| \ot_{|fF_n|}|A| \;$;

\item  the space $C(M;A)$ has the weak topology with respect
to the filtration  $C_k(M;A) = Im( \coprod_{i \leq k} fC(M)_i
\times_{\s_i} A^i) ,\;$
$k \in \N$.
\end{enumerate}

\end{prop}
We give some examples of configuration spaces with summable labels.
Let us denote by $\hat{A}^n$ the completion of a partial $n$-monoid $A$.
\begin{prop} \label{euclidean}
If $A$ is a partial $n$-monoid, then
there is a strong deformation retraction $w_A:C(\R^n;A) \to \hat{A}^n$.
\end{prop}
\begin{proof}
It is sufficient to observe that there is a deformation retraction
of right $F_n$-modules $w:C(\R^n) \to F_n$. If an element $x \in C(\R^n;A)$
is represented by a finite number of labelled trees $\tau_1,\dots,\tau_k$
based at distinct
points $P_1,\dots,P_k \in \R^n$, then $w_A(x)$ is represented by
the single tree obtained by connecting $\tau_1,\dots,\tau_k$ to a root
vertex labelled by the class $[P_1,\dots,P_k] \in C^0_k(\R^n)/G(n) \subset
F_n(k)$.

\end{proof}

\begin{ex}
If $M$ is a discrete partial monoid, then
$\hat{M}^1$ has the homotopy type of its monoid completion.
If $M$ is abelian then $\hat{M}^\infty$ has the homotopy type of its 
abelian monoid completion.
\end{ex}

\begin{defi}
Let $A$ be a partial abelian monoid and $M$ an $n$-manifold.
We denote by $C^0(M;A)$ the quotient of $\coprod_k C^0_k(M) \times_{\s_k}
A^k$ under the following relation $\sim$:  if $(m_1,\dots,m_k) \in
C^0_k(M)$, $a_1,\dots,a_k \in A$,
$m_1 = m_2$ and $a_1 + a_2$ is defined, then
$$(m_1,\dots,m_k;a_1,\dots,a_k) \sim (m_2,\dots,m_k;a_1+a_2,\dots,a_k).$$
\end{defi}

\begin{lem} \label{tatu}
If $A$ is a partial abelian monoid, and $M$ is an $n$-dimensional
open
manifold, then the inclusion
$C^0(M;A) \hook C(M;A)$ is a weak equivalence.
\end{lem}

\begin{proof}

The proof makes use of the fact that a copy of the manifold
with corners $C_k(M)$ lies inside its interior $C^0_k(M)$, so the retraction
$r: C_k(M) \to C^0_k(M)$ is a $\s_k$-equivariant homeomorphism onto its
image.
We compare via this retraction the pushout diagram for $C^0_k(M;A)$

$$\begin{diagram}
C_k^0(M) \times_{\s_k} T_k(X,x_0) \cup (C^0(M) \times_{\tau} Comp(A))_k
 & \rTo & C^0_{k-1}(M;A) \\
\dHook     &             & \dHook        \\
C^0_k(M) \times_{\s_k} A^k  \cup (C^0(M) \times_{\tau} Comp(A))_k &
 \rTo    &      C^0_k(M;A) ,\\
\end{diagram} $$

and the pushout diagram for $C_k(M;A)$

$$\begin{diagram}
C_k(M) \times_{\s_k} T_k(X,x_0) \cup (C^0(M) \times_{\tau} Comp(A))_k &
\rTo & C_{k-1}(M;A) \\
\dHook     &             & \dHook         \\
  C_k(M) \times_{\s_k} A^k    &     \rTo    &      C_k(M;A) .\\
\end{diagram} $$

Here we denote by $(C^0(M) \times_\tau Comp(A))_k$
 the subspace of $C(M)_k \times_{\s_k} A^k$ of those labelled configurations
such that several points are concentrated in the same macroscopic location
if and only if their labels are summable.
The inclusion of the space on the left hand top corner of the first diagram
into that of the second diagram is a homotopy equivalence, because
$r$ induces a common retraction onto a copy of the second space.
The same holds for the spaces on the left hand bottom corner.
We conclude by induction and the gluing lemma \cite{Brown}.

\end{proof}

If we regard a pointed space $(A,a_0)$ as a partial
abelian monoid with $x+a_0 = x$ as the only defined sums, for $x \in A$,
then $C^0(M;A)$ is the configuration space with labels studied in
\cite{Bodig}.

\begin{cor} \label{stanco}
Let $(A,a_0)$ be a well-pointed space.
Then for any open $n$-manifold $M$ there is a weak equivalence
 $C^0(M;A) \simeq C(M;A)$.
\end{cor}
\begin{proof}
The space $A$ is a partial $n$-monoid
by $Comp = \coprod_k F_n(k) \times_{\s_k} \vee_{i=1}^k   A $.
\end{proof}
\

For some background about toric varieties we refer to \cite{Fulton}.

\begin{cor} \label{Guest}
If $V$ is a projective
toric variety such that $H_2(V)$ is torsion free, then there exists a partial discrete abelian
monoid $\Delta_V$, such that the union of some components of
$\hat{(\Delta_V)}^2$ is homotopy equivalent
to the space $Rat(V)$
of based rational curves on $V$.
\end{cor}
\begin{proof}
Guest has shown in \cite{Guest} that if $\Delta_V$ is the fan associated
to the variety $V$ \cite{Fulton} then the union of some components
of $C^0(\R^2;\Delta_V)$ is homeomorphic to $Rat(V)$. The corollary follows
from the theorem and from proposition \ref{euclidean}.
\end{proof}

\begin{rem}
It is possible to define labelled configurations with support in a manifold with
corners $M$. It is sufficient to choose an embedding $M \hook M'$, with
$M'$ open, consider the right $F_n$-submodule
$C(M) \hook C(M')$ of configurations macroscopically
located at points of $M$, and carry through the discussion
as for open manifolds. 
\end{rem}

\section{The relative case}
We define relative labelled configuration spaces on relative manifolds.

Let $(X,x_0)$ be a pointed topological space.
Let $M$ be a manifold with corners and $N \hook M$ a cofibration such that
$M - N$ is an open manifold.
We obtain easily from \ref{Kont} that
each element $c \in
C(M;X)$ is uniquely determined by
 a finite set
$S(c) \subset M$, and for each $P \in S(c)$ a labelled tree
$T_P$ as in \ref{Kont},
with the only difference that the twigs of the tree are labelled by
$X-x_0$.

\begin{defi}
The based space $C(M,N)(X)$ is
the quotient  $C(M;X)/ \sim$ by the equivalence
relation such that $a \sim a'$
if and only if 
$S(a) \cap (M-N) = S(a') \cap (M-N)$ and the trees indexed by these
intersections coincide. The base point is the class $[a]$ such that
$S(a) \subset N$.
\end{defi}

If we regard pointed spaces as partial $n$-monoids, then the
$n$-monoid completion induces a monad
$(F_n^*,\eta_*,\mu_*)$  on the 
category of pointed compactly generated
spaces $\mathcal{CG}_*$.
Each element in the completion $F_n^*(X) = \hat{X}^n$ is represented
by a tree with vertex labels in $F_n$ and twigs labels in $X-x_0$.
The product $\mu_*$ is given by grafting of trees, and the unit $\eta_*$
sends an element $x$ to the trivial tree labelled by $x$.

\begin{prop}
If $M$ is a parallelizable $n$-manifold, and $N \hook M$ is a cofibration
such that $M-N$ is open,
then the functor $C(M,N)$ has a structure of right algebra over
$F_n^*$.
\end{prop}

\begin{proof}
We need to exhibit a natural transformation
$\lambda: C(M,N) F_n^* \to C(M,N)$ such that
$\lambda \circ C(M,N)\eta_*$ is the identity and
the diagram
$$\begin{diagram}
C(M,N) F_n^* F_n^* & \rTo^{C(M,N) \mu_*} & C(M,N) F_n^*  \\
\dTo_{\lambda F_n^*} &                      & \dTo^{\lambda} \\
C(M,N) F_n^*            & \rTo^{\lambda}      & C(M,N)         .\\
\end{diagram} $$
commutes. The morphism $\lambda$ is obtained by grafting of trees.
\end{proof}

\begin{defi}
If $(A,\rho)$ is an $n$-monoid, and $M,N$ are as before,
then the space $C(M,N;A)$ of configurations
in $(M,N)$ with summable labels
in $A$   is the coequalizer
$$C(M,N) F_n^*( A)
 \pile{\rTo^{C(M,N) \rho} \\
\rTo_{\lambda A} }
C(M,N) A \rDotsto C(M,N;A) \; .$$
\end{defi}

\begin{defi} A partial $n$-monoid $A$ is good if
the inclusion $Comp(A) \to F_n(A)$ is a cofibration, and the partial
composition $\rho: Comp(A) \to A$ induces a map on the quotient
$Comp^*(A) \subset F_n^*(A)$ of $Comp(A)$.
\end{defi}
The definition of a good framed partial $n$-monoid is similar.
From now on we will assume implicitly that all partial (framed)
$n$-monoids are good.


By means of the framed Fulton-MacPherson operad we can define similarly
$C(M,N;A)$, if $M$ is an $n$-dimensional
manifold with corners,  $N \hook M$ is a cofibration, and $A$ is
a good partial framed $n$-monoid, and as in \ref{pish} we obtain:

\begin{prop}
Define a filtration so that $[a] \in C_k(M,N;A)$ if and only if $k$ is
the number of twigs of trees in 
$S(a) \cap (M-N)$. Then $C(M,N;A)$ has the weak topology
with respect to the filtration and it is compactly generated.
\end{prop}

\begin{defi}
If $A$ is a partial framed $n$-monoid, then 
$B_k(A) = C( (I^k,\partial I^k) \times I^{n-k} ;A)$
for $i=1 ,\dots, n$.
\end{defi} 

If $A$ is a partial abelian monoid and $(M,N)$ is any pair then
we define the relative labelled configuration space
$C^0(M,N;A)$ as quotient of $C^0(M;A)$, by identifying configurations
that coincide on $M-N$.
We state the relative version of \ref{tatu}.
\begin{prop}
If $M$ is a manifold, $N \hook M$ is a cofibration, and $M-N$ is open, then
there
is a weak equivalence $C^0(M,N;A) \simeq C(M,N;A)$.
\end{prop}

\begin{proof}
The proof is similar to that of \ref{tatu}.
In this case we use for each $k$ a $\s_k$-equivariant retraction
$r_k: C_k(M) \to C_k^0(M)$ such that $r_k$ preserves
$b^{-1} (\overline{M^k - N^k  }    )$, where $b: C_k(M) \to M^k$ is the
blowdown.
\end{proof}

\begin{cor} \label{b2}
If $V$ is a projective toric variety such that
$H_2(V)$ is torsion free,
with torus $T$ and fan $\Delta_V$, then
there is a weak equivalence $B_2(\Delta_V) \simeq V \times_T ET$.
\end{cor}
\begin{proof}
Guest has shown in \cite{Guest} that $V \times_T ET$ is homotopy equivalent
to

$C^0(I^2,\de I^2; \Delta_V)$.
\end{proof}
The relative version of \ref{stanco} is:

\begin{cor} \label{dosta}
For any well-pointed space $X$ there is a weak equivalence 

$C^0(M,N;X) \simeq C(M,N;X)$.
\end{cor}

\begin{cor} \label{corolla}
Let $X$ be a well-pointed space considered as partial $n$-monoid.
Then there is a weak equivalence
$\s^n(X) \rHook^\simeq B_n(X)$.
\end{cor}
\begin{proof}
The space of open configurations $C^0(I^n,\de I^n;X)$ 
retracts onto $\s^n(X)$, considered as space of configurations 
of a single labelled point in $(I^n,\de I^n)$.
The retraction is achieved  \cite{dusa} by pushing radially
the particles away onto the boundary.
But the inclusion $C^0(I^n,\de I^n;X) \hook C(I^n,\de I^n;X) = B_n(X)$
is a weak equivalence by \ref{dosta}.
\end{proof}

By means of configuration spaces
we obtain the classifying space constructed by Stasheff.
\begin{prop}
Let $(A,a_0)$ be a well-pointed $A_\infty$ space.
The quotient map
$C(I,\{0\};A) \to C(I,\partial I;A) = B_1(A)$ is
canonically homeomorphic to the universal arrow $E(A) \to B(A)$.
\end{prop}

\begin{proof}
It is sufficient to carry out the discussion in the non-symmetric
case: in fact 
$C_k([0,1])=S_k([0,1]) \times \s_k$, where $S_k([0,1])$
compactifies the
space
of strictly ordered maps from $\{1,\dots,k\}$ to $I=[0,1]$.

Let $S_k(I)\{0,1\} \subseteq S_k(I)$ be 
the closure of the subspace of maps $\alpha: \{1,\dots,k\} \to I$ 
such that $\alpha(1)=0$, $\alpha(k)=1$.
Its elements are described by appropriate trees as in
\ref{Kont}.
For $k \geq 2$, we have homeomorphisms
 $r: S_k(I)\{0,1\} \leftrightarrows S_k(0): j$, where
$S_k(0)$ is
  the space of
configurations in $\R$ macroscopically concentrated at the point
$0$. 

If $\alpha_i \to \alpha \in S_k(I)\{0,1\}$, $\alpha_i \in 
C^0_k(I)$, then $r(\alpha)=\lim_i \frac{\alpha_i
-\alpha_i(0)}
{i(\alpha_i(1)-\alpha_i(0))} $.

If $\beta_i \to \beta \in S_k(0)$, $\beta_i \in C^0_k(\R)$, then
$j(\beta) = \lim_i \frac{\beta_i - \beta_i(0)}{\beta_i(1)-\beta_i(0)}$.

We have seen in \ref{Kont} that $S_k(0) = K_k$ is the associahedron.
Under the identification $K_k \cong S_k(I)\{0,1\}$
the Stasheff space $B(A)$ is defined to be the
quotient of $\coprod S_k(I)\{0,1\} \times A^{k-2}$, seen as space
of forests labelled by $A$,
under the following steps:
\begin{enumerate}
\item We replace a tree on $i$ twigs by a point having as label the
action of the tree on its twigs via $K_i \times A^i \to A$.
\item We can cut twigs labelled  by $a_0$.
\item We identify any two labelled forests coinciding
outside 0 and 1.
\end{enumerate}
But this quotient is exactly
$B_1(A) = C(I,\partial I;A)$.
In a similar way one shows that $E(A) = \coprod S_{k}(I)\{0,1\} \times
A^{k-1} / \sim$
is homeomorphic to  $C(I,\{0\};A)$. In this case in 3
we identify forests coinciding outside 0.
\end{proof}

\section{Approximation theorems}

We say that a partial framed $n$-monoid $A$ has homotopy inverse if
the H-space $\hat{A}^n$ has homotopy inverse.

\begin{lem} \label{isotopy}
Let $M$ be a  connected compact  $n$-manifold,
$M' \subset M$ a compact $n$-submanifold,
$N \subset M$ a closed submanifold, 
and $A$ a partial framed $n$-monoid.
Suppose that either $A$ has a homotopy inverse
or the pair $(M',N\cap M')$ is connected.
Then  there is
a quasifibration

$C(M',N\cap M';A) \rTo C(M,N;A) \rTo^\pi C(M,M' \cup N;A)$.
\end{lem}
This holds in particular if $A$ is path connected.
\begin{proof}
We follow the proof of proposition 3.1 in  \cite{dusa}.
The space $C(M,M' \cup N;A)$ has a filtration by 
$C_k := C_k(M,M' \cup N;A)$.
There is a homeomorphism 
$\alpha_k: \pi^{-1} (C_k - C_{k-1})  \cong
C(M',N \cap M';A) \times (C_k - C_{k-1})$ such that
$\pi \alpha_k^{-1}$ is the projection onto the factor $C_k - C_{k-1}$. 
Choose a collared neighbourhood $U$ of $M'$ in $M$
and a smooth isotopy retraction
$r:U \to M'$ such that
$r(U \cap N) \subset N$.
For each $k$ there is an open neighbourhood $U_k$ of $C_k$ in $C_{k+1}$
such that $r$ induces a smooth isotopy retraction
$r_k:U_k \times I\to C_k$, and a smooth isotopy retraction
$\tilde{r}_k: \pi^{-1}(U_k) \times I \to \pi^{-1}(C_k) $ covering $r_k$.
For any point $P \in U_k$ we need to show that the restriction
$t: \pi^{-1}(P)  \to \pi^{-1}(r_1(P))$ of $\tilde{r}_1$ is a weak homotopy
equivalence. If we identify domain and range of $t$ to $C(M',N \cap M';A)$ by
$\alpha_k$, then  $t$ 
pushes the labelled particles away from $N$, and adds a finite
set of trees $T$ in proximity to $N$. But if the pair $(M',N \cap M')$ is
connected, then the trees in $T$ can be moved continuously to $N$, where
they vanish, and $t$ is homotopic to a homeomorphism.
On the other hand, if $A$ has a homotopy inverse, then $t$ has
a homotopy inverse
that pushes the particles away from $N$ and adds some homotopy inverses
of the trees in $T$ in proximity to $N$.
\end{proof}

\begin{prop} \label{mai}
Let $A$ be a partial framed $n$-monoid.
Then for $i =1, \dots, n$ there are maps
$s_i: B_{i-1}(A) \rTo \Omega B_i(A)$, such that
$s_i$ is a weak homotopy equivalence for $i>1$,  and
$s_1$ is a weak homotopy equivalence
if $A$ has a homotopy inverse.
\end{prop}

\begin{proof}
Note that $B_0(A)$ is homotopic to the framed $n$-monoid completion of $A$.
For each $i$ the base point of $B_i(A)$ is the empty configuration.
The translation $\tau_1(t): I^n \to \R \times I^{n-1}$
of the first coordinate
by $t$ induces a map $\pi_{\tau_1(t)}: B_0(A) \to B_1(A)$, composite of
the induced map
$C(I^n;A) \rTo^{C(\tau_i(t);A)} C(\R \times I^{n-1};A)$ and the projection
$C(\R \times I^{n-1};A) \to C((I,\de I) \times I^{n-1};A)$.
 Then the `scanning' map $s_1$
is defined for $x \in B_0(A)=C(I^n;A)$
by $s_1(x)(t)= \pi_{ \tau_1(2t-1)}(x) \in B_1(A) $.
For $i>0$ 
the translation
of the $(i+1)$-th coordinate by $t$ induces similarly
a map $\pi_{\tau_{i+1}(t)}: B_{i}(A) \to B_{i+1}(A)$, and 
$s_{i+1}:B_i(A) \to \Omega B_{i+1}(A)$ is given by
$s_{i+1}(x)(t) = \pi_{\tau_{i+1}(2t-1)} (x)$.
We define $M= I^k \times [0,2] \times I^{n-k-1}$, $N=(\de I^k \times
[0,2] \times I^{n-k-1}) \cup  (I^k \times 0 \times I^{n-k-1})$, and
we identify $B_k(A)$ to
$C(I^k \times [1,2] \times I^{n-k-1}, \de I^k \times [1,2] \times
I^{n-k-1};A)$ via $\tau_{k+1}(1)$.
We consider
for $1 \leq k \leq n-1$ 
a
commutative diagram

$$\begin{diagram} B_k(A) & \rTo &  C(M,N;A) & \rTo &  B_{k+1}(A) \\ 
           \dTo^{s_{k+1}} &     &   \dTo^s  &        &   \dEqual \\
 \Omega B_{k+1}(A) & \rTo & PB_{k+1}(A) & \to & B_{k+1}(A)\;.
 \end{diagram} $$
The top row is a quasifibration and the bottom row a fibration.
The scanning map $s$ is defined
on the total space $C(M,N;A)$ by $s(x)(t) = \pi_{\tau_{k+1}(2t)}(x)$
and is consistent with $s_{k+1}$.
Now the space $C(M,N;A)$ is contractible. In fact by excision 
$C(M,N;A) \cong C(M',N';A)$, with $M'=\R^k \times (-\infty,2] \times
\R^{n-k-1}$ and $N' = M' -(M-N)$. Moreover there is a smooth
isotopy $H_t:(M,N) \to (M',N')$, such that
$H_0$ is the inclusion and $H_1(M) \subset N'$. For example define
$H_t$ as the dilatation by $3t$ centered in $(\frac{1}{2},\dots,
\frac{1}{2},3,\frac{1}{2}, \dots, \frac{1}{2})$, with $3$ at the 
$(k+1)$-st position. We conclude by comparing the long exact sequences
in homotopy and by induction on $k$.
\end{proof}

The spaces
$B_0(A)=C(I^n;A)$
and $\Omega^n B_n(A)$ are both $fD_n$-algebras. 
The map $s: B_0(A) \to \Omega^n B_n(A)$
constructed by
looping and composing the
scanning maps in
proposition \ref{mai} can be extended to
a homotopy $fD_n$-morphism 
by rescaling
suitably the scanning maps on the labels of trees in $U_{fD_n}(B_0(A))$.
By \ref{mai} we obtain:

\begin{thm} 
If $A$ is a partial framed $n$-monoid, then 
$s: B_0(A) \to \Omega^n B_n(A)$ is the group completion.
If $A$
has homotopy inverse, then $s$ is a weak homotopy equivalence.
\end{thm}

Actually $s$ is the group completion in the
homotopy category of $fD_n$-algebras.

\begin{cor} {\rm \cite{May}}
If $X$ is a well-pointed space, then 
 $s: C^0(\R^n;X) \to 
\Omega^n \Sigma^n X$ is the group completion. If $X$ is path connected,
then $s$ is a weak homotopy equivalence.
\end{cor}
\begin{proof}
Consider $X$ as a partial $n$-monoid as in corollary \ref{stanco}.
Now $B_0(A) = C(I^n;A) \simeq \hat{A}^n $ by the same argument
of proposition \ref{euclidean} . Moreover $\hat{A}^n \simeq C(\R^n;A)$ by
\ref{euclidean},
$C(\R^n;A) \simeq C^0(\R^n;A)$ by \ref{stanco} and $\s^n X \simeq B_n(X)$ by \ref{corolla}.
Now we can apply the theorem.
\end{proof}
\begin{cor} {\rm \cite{Guest}}
If $V$ is a projective toric variety such that $H_2(V)$ is torsion free,
then $s: Rat(V) \to \Omega^2(V)$
is the group completion.
\end{cor}
\begin{proof}
Apply corollaries \ref{Guest} and \ref{b2},
and restrict to the relevant components.
\end{proof}

Given a manifold $M$, and its tangent bundle $\tau$, there is
a bundle $\gamma = C(\tau,\partial \tau;A)$ on $M$ with fiber
$B_n(A) = C(I^n,\partial I^n;A)$, consisting of relative 
fiberwise configurations in the fiberwise one-point compactification
modulo the section at infinity
$(\hat{\tau},\infty)$.
Whether $\de M$ is empty or not we can define a map 
$s: C(M,\de M;A) \to \Gamma(M;B_n A)$ to the space of sections
of $\gamma$.
Note that if $M$ is parallelizable
then $\Gamma(M;B_nA) = Map(M;B_nA)$.
The scanning map $s$ is constructed by the
exponential map: if $x \in C(M,\de M;A)$, then $s(x)$ sends
a point $P \in M$
to the restriction of $x$ to a small disc neighbourhood of $P$
modulo its boundary.

\begin{thm}
Let $A$ be a partial framed $n$-monoid.
Let $M$ be a compact connected $n$-manifold with boundary.
Then the scanning map $s: C(M,\de M;A) \to \Gamma(M;B_n A)$
is a weak homotopy equivalence.
If $A$ has homotopy inverse and $N$  is a compact connected $n$-manifold without boundary
then $s:C(N;A) \to \Gamma(N;B_nA)$ is a weak homotopy equivalence.
\end{thm}

\begin{proof}
We follow the proof of 10.4 in \cite{Kallel}.
There is a finite handle decomposition of $M$ with no handles of index $n$.
If $M'$ is obtained from $M''$ by attaching a handle $H$ of index $i$
, then we apply lemma \ref{isotopy}
and we obtain a quasifibration
$C(H, \overline{ \de H - \de H \cap M''} ;A)    \to C(M', \de M';A) \to C(M'',\de M'';A)$.
On the other hand  we have a fibration
$\Gamma(H/(H \cap M'');B_nA)  \to \Gamma(M';B_nA)
\to \Gamma(M'';B_nA)$.
But  $C(H, \overline{ \de H - \de H \cap M''};A ) \cong B_{n-i}(A)$, and
$\Gamma(H/(H \cap M'');B_nA) \simeq \Omega^i  B_n(A)$.
We compare the two sequences by the scanning maps and we conclude by
proposition \ref{mai} and
induction on the number of handles.
In the case of $N$ we have even a handle of index $n$ and we
apply the second part of proposition \ref{mai}.
\end{proof}

\begin{cor}
If $X$ is a well-pointed path connected space then
$s: C(S^1;X) \to Map(S^1, \Sigma X)$ is a 
weak homotopy equivalence.
\end{cor}

\begin{proof}
We consider $X$ as partial framed $1$-monoid as in corollary \ref{stanco}.
By corollary \ref{corolla} $B_1(X) \simeq \s X$. We apply the
second part of the theorem, and note that
$\Gamma(S^1,\s X) \simeq Map(S^1,\s X)$
because $S^1$ is parallelizable.
\end{proof}

This answers a question raised by Stasheff in \cite{operadchik}  p. 10.
The analogous result for $C^0(S^1;X)$ is in \cite{Bodig}.

Any partial framed $n$-monoid gives an approximation theorem for
mapping spaces, and the homotopy theorist is tempted to discover new examples. 
It might be worth considering colimits of abelian monoids in the category
of $n$-monoids.

\

Mathematisches Institut

Universit\"at Bonn

Beringstrasse 1

53115 Bonn

Germany
 
e-mail: salvator@math.uni-bonn.de

\end{document}